\documentclass[a4paper]{article}%
\usepackage{graphicx}
\usepackage{geometry}
\usepackage{amsmath}
\usepackage{amssymb}
\usepackage{amsfonts}
\usepackage{amsthm,amscd}%
\providecommand{\U}[1]{\protect\rule{.1in}{.1in}}
\def\figurename{Figure}
\makeatletter
\renewcommand{\fnum@figure}[1]{\figurename~\thefigure.}
\makeatother
\def\tablename{Table}
\makeatletter
\renewcommand{\fnum@table}[1]{\tablename~\thetable.}
\makeatother
\def \bop {\noindent\textbf{Proof. }}
\def \eop {\hbox{}\nobreak\hfill
\vrule width 2mm height 2mm depth 0mm
\par \goodbreak \smallskip}
\newtheorem{theorem}{Theorem}[section]
\newtheorem{lemma}[theorem]{Lemma}
\newtheorem{corollary}[theorem]{Corollary}
\newtheorem{proposition}[theorem]{Proposition}
\theoremstyle{definition}
\newtheorem{definition}[theorem]{Definition}

\theoremstyle{remark}
\newtheorem{remark}[theorem]{Remark}
\numberwithin{equation}{section}

\begin{document}

\title{On the relaxed mean-field stochastic control problem \thanks{Partially
supported by the French Algerian Cooperation Program, Tassili 13\ MDU 887. }}
\author{{Khaled Bahlali}\thanks{ Laboratoire IMATH, Universit\'{e} du Sud-Toulon-Var,
B.P 20132, 83957 La Garde Cedex 05, France. E-mail: bahlali@univ-tln.fr}
\and {Meriem Mezerdi }\thanks{ Laboratory of Applied Mathematics, University of
Biskra, Po. Box 145, Biskra (07000), Algeria. E-mail: m\_mezerdi@yahoo.fr}
\and Brahim Mezerdi\thanks{ Laboratory of Applied Mathematics, University of
Biskra, Po. Box 145, Biskra (07000), Algeria. E-mail: bmezerdi@yahoo.fr}}
\maketitle

\begin{abstract}
This paper is concerned with optimal control problems for systems governed by
mean-field stochastic differential equation, in which the control enters both
the drift and the diffusion coefficient. We prove that the relaxed state
process, associated with measure valued controls, is governed by an orthogonal
martingale measure rather that a Brownian motion. In particular, we show by a
counter example that replacing the drift and diffusion coefficient by their
relaxed counterparts does not define a true relaxed control problem. We
establish the existence of an optimal relaxed control, which can be
approximated by a sequence of strict controls. Moreover under some convexity
conditions, we show that the optimal control is realized by a strict control.

\end{abstract}



\section{Introduction}

In this paper, we deal with optimal control of systems driven by mean-field
stochastic differential equations (MFSDE) of the form%

\[
\left\{
\begin{array}
[c]{l}%
dX_{t}=b(t,X_{t},E(\Psi(X_{t})),u_{t})dt+\sigma(t,X_{t},E(\Phi(X_{t}%
)),u_{t})dW_{t}\\
X_{0}=x.
\end{array}
\right.
\]

MFSDEs are obtained as limits of some interacting particle systems. This kind
of approximation result is called "\textit{propagation of chaos}", which says
that when the number of particles (players) tends to infinity, the equations
defining the evolution of the particles could be replaced by a single
equation, called the McKean-Vlasov equation. This mean-field equation,
represents in some sense the average behavior of the infinite number of
particles, see \cite{Sn, JMW} for details. Since the earlier papers
$\cite{LasLio}$, $\cite{HMC}$, mean-field control theory has raised a lot of
interest, motivated by applications to various fields such as game theory,
mathematical finance, communications networks, management of oil
ressources.\ The typical example is the continuous-time Markowitz's
mean-variance portfolio selection model, where one should minimize an
objective function involving a quadratic function of the expectation, due to
the variance term. The main drawback, when dealing with such mean-field
stochastic control problems, is that the Bellmann principle of optimality does
not hold. For this kind of problems, the stochastic maximum principle,
provides a powerful tool to solve them, see \cite{AndDje, BucDjeLi, CarDel,
ChiMez, ELN, Li, Temb}. One can refer also to the recent book \cite{Bens} and
the references therein.

Our main goal in this work is to investigate existence of optimal
controls.\ As it is well known for classical control problems, in the absence
of Fillipov convexity conditions, an optimal strict control may fail to exist.
In this case, the set of strict controls should be embedded into a wider class
of measure valued controls, called relaxed controls. This class enjoys good
compactness and convexity properties. The problem now is to define precisely
the MFSDE associated to a relaxed control. At first look, one is tempted as in
\cite{Cha}, to replace simply the drift and diffusion coefficient by their
relaxed couterparts ie: the integrals of the drift and diffusion coefficient
with respect to the relaxed control, adopting the same method as in
deterministic control. As it will be shown in a simple counter example, the
suggested "relaxed" state equation is not continuous with respect to the
control variable. This implies in particular that the value functions for the
original and relaxed problems are not the same. In addition, there is no mean
to prove the existence of an optimal control for this model. So that the
proposed "relaxed" model in \cite{Cha} is not a true extension of the original
control problem. The fundamental reason is that one has to relax the quadratic
variation, of the stochastic integral part of the state equation, which is a
Lebesgue integral, rather than the stochastic integral itself. Roughly
speaking, the idea is to relax the generator of the process, which is
intimately linked to the weak solutions of the relaxed stochastic equation,
rather than the equation itself. As it will be shown, the stochastic equation
associated with the relaxed generator will be governed by a continuous
orthogonal martingale measure, rather than a Brownian motion. For this model,
we prove that the strict and relaxed control problems have the same value
function and that an optimal relaxed control exists. Our result extends in
particular \cite{BahDjeMez, BDM, EKNJ, BahMez1, BahMez} to mean field controls
and \cite{BMM1} to the case of a MFSDE with a controlled diffusion
coefficient. The proof is based on tightness properties of the underlying
processes and Skorokhod selection theorem. Moreover, due to the compactness of
the action space, we show that the relaxed control could be choosen among the
so-called sliding controls, which are convex combinations of Dirac measures.
As a consequence and under the so-called Fillipov convexity condition, the
optimal relaxed control is shown to be strict.

\section{Existence of optimal relaxed controls}

\subsection{Controlled mean field stochastic differential equations}

Let $\left(  W_{t}\right)  $ is a $d$-dimensional Brownian motion, defined on
a probability space $(\Omega,\mathcal{F},P),$ endowed with a filtration
$\left(  \mathcal{F}_{t}\right)  ,$ satisfying the usual conditions. Let
$\mathbb{A}$ be some compact metric space called the control set.

We study the existence of optimal controls for systems driven non linear
mean-field stochastic differential equations (MFSDE in short) of the form:%

\begin{equation}
\left\{
\begin{array}
[c]{l}%
dX_{t}=b(t,X_{t},E(\Psi(X_{t})),u_{t})dt+\sigma(t,X_{t},E(\Phi(X_{t}%
)),u_{t})dW_{t}\\
X_{0}=x.
\end{array}
\right.  \label{MFSDE1}%
\end{equation}
\bigskip

The cost functional over the time interval $\left[  0,T\right]  $ is given by%

\begin{equation}
J(u)=E\left(
{\displaystyle\int\limits_{0}^{T}}
h(t,X_{t},E(\varphi(X_{t})),u_{t}\right)  dt+g(X_{T},E(\lambda(X_{T}))),
\label{COST1}%
\end{equation}
where $b,$ $\sigma,$ $l,$ $h$, $g$ and $\psi$ are given functions. The control
variable $u_{t}$ called a strict control, is a measurable, $\mathcal{F}_{t}-$
adapted process with values in the action space $\mathbb{A}$. We denote
$\mathcal{U}_{ad}$ the space of strict controls. Let us point out that the
probability space and Brownian motion may change with the control $u.$

The objective is to minimize the cost functional $J(u)$ over the space
$\mathcal{U}_{ad}$ ie: find $u^{\ast}\in$ $\mathcal{U}_{ad}$ such that
$J(u^{\ast})=\inf\left\{  J(u);u\in\mathcal{U}_{ad}\right\}  .$

The following assumptions will be in force throughout this paper.

$\mathbf{(H}_{\mathbf{1}}\mathbf{)}$ Assume that%

\begin{equation}%
\begin{array}
[c]{c}%
b:\left[  0,T\right]  \times\mathbb{R}^{d}\times\mathbb{R}^{d}\times
\mathbb{A}\longrightarrow\mathbb{R}^{d}\\
\sigma:\left[  0,T\right]  \times\mathbb{R}^{d}\times\mathbb{R}^{d}%
\times\mathbb{A}\longrightarrow\mathbb{R}^{d}\otimes\mathbb{R}^{d}\\
\Psi:\mathbb{R}^{d}\longrightarrow\mathbb{R}^{d},\Phi:\mathbb{R}%
^{d}\longrightarrow\mathbb{R}^{d}%
\end{array}
\end{equation}

are bounded continuous functions such that $b(t,.,.,a),$ $\sigma
(t,.,.,a),\Psi(.)$ and $\Phi(.)$ are Lipschitz continuous, uniformly in
$(t,a)$.

$\mathbf{(H}_{\mathbf{2}}\mathbf{)}$ Assume that%

\begin{equation}%
\begin{array}
[c]{l}%
h:\left[  0,T\right]  \times\mathbb{R}^{d}\times\mathbb{R}^{d}\times
\mathbb{A}\longrightarrow\mathbb{R}\\
g:\mathbb{R}^{d}\times\mathbb{R}^{d}\longrightarrow\mathbb{R}\\
\varphi:\mathbb{R}^{d}\longrightarrow\mathbb{R}^{d}\\
\lambda:\mathbb{R}^{d}\longrightarrow\mathbb{R}^{d}%
\end{array}
\end{equation}

are bounded continuous functions and $h(t,.,.,a)$ is Lipschiz continuous
uniformly in $(t,a)$.

\begin{proposition}
Under assumption $\mathbf{(H}_{\mathbf{1}}\mathbf{)}$ the MFSDE (\ref{MFSDE1})
has a unique strong solution, such that for every $p>0$ we have $E(\left\vert
X_{t}\right\vert ^{p})<+\infty.$
\end{proposition}

\bop

Let us define $\overline{b}(t,x,\mu,a)$ and $\overline{\sigma}(t,x,\mu,a)$ on
$\left[  0,T\right]  \times\mathbb{R}^{d}\times\mathbb{P}_{2}(\mathbb{R}%
^{d})\times\mathbb{A}$ by

$\overline{b}(t,x,\mu,a)=b(t,x,%
{\displaystyle\int}
\Psi(x)d\mu(x),a)$ and $\overline{\sigma}(t,x,\mu,a)=\sigma(t,x,%
{\displaystyle\int}
\Phi(x)d\mu(x),a),$ where $\mathbb{P}_{2}(\mathbb{R}^{d})$ denotes the space
of probability measures in $\mathbb{R}^{d},$ having a finite second order
moment$.$

According to \cite{JMW} Prop.1.2, it is sufficient to check that $\overline
{b}$ and $\overline{\sigma}$ are Lipschitz in $\left(  x,\mu\right)  $ where
$\mathbb{P}_{2}(\mathbb{R}^{d})$ is equipped with the 2-Wasserstein metric
$d\left(  \mu,\nu\right)  =\inf\left\{  \left(  E^{Q}\left\vert X-Y\right\vert
^{2}\right)  ^{1/2};Q\in\mathbb{P}_{2}(\mathbb{R}^{d}\times\mathbb{R}%
^{d}),\text{ with marginals }\mu,\nu\right\}  .$ This follows from the
Lipschitz continuity of $b$ and $\sigma$ with respect to $(x,y).$

Using similar techniques as in \cite{JMW} Prop.1.2, it holds that for each
$p>0$, $E(\left\vert X_{t}\right\vert ^{p})<+\infty.$ \eop

\bigskip

\textbf{Examples of mean-field control problems}

1) \textbf{Example 1:} \textbf{The mean-variance portfolio selection problem}

Consider a financial market, in which two securities are traded continuously.
The first is a bond, with price $S_{t}^{0}$ at time $t\in\left[  0,T\right]  $
governed by
\[
dS_{t}^{0}=S_{t}^{0}\rho_{t}dt,\text{ }S_{0}^{0}=s_{0}>0.
\]

The second is a stock with unit price $S_{t}^{1}$ at time $t\in\left[
0,T\right]  $ governed by%
\[
dS_{t}^{1}=S_{t-}^{1}\left(  b_{t}dt+\sigma_{t}dB_{t}\right)  ,\text{ }%
S_{0}^{1}=s^{1}>0.
\]

The coefficients $\rho_{t}>0,b_{t},\sigma_{t}$ are deterministic and bounded
functions. For an investor, a portfolio $\pi$ is a process representing the
amount of money invested in the stock$.$ The wealth process $x^{x_{0},\pi}$
corresponding to initial capital $x_{0}>0,$ and portfolio $\pi$, satisfies
then the equation%
\[
\left\{
\begin{array}
[c]{l}%
dx_{t}=\left(  \rho_{t}x_{t}+\pi_{t}\left(  b_{t}-\rho_{t}\right)  \right)
dt+\pi_{t}\sigma_{t}dB_{t},\text{ for }t\in\left[  0,T\right]  ,\\
x_{0}=x.
\end{array}
\right.
\]

The objective is to maximize the mean terminal wealth $\mathbb{E}\left[
x_{T}^{\pi}\right]  ,$ and at the same time, to minimize the variance of the
terminal wealth $\mathrm{Var}\left[  x_{T}^{\pi}\right]  ,$ over controls
$\pi$ valued in $%
\mathbb{R}
$. Then, the mean-variance portfolio optimization problem is to minimize the
cost $J$, given by
\[
J\left(  \pi\right)  =-\mathbb{E}\left[  x_{T}\right]  +\mu\mathrm{Var}\left[
x_{T}\right]  ,
\]
subject to $\left(  4.3\right)  $, where $\mu>0.$ The admissible portfolio is
assumed to be progressively measurable square integrable process$,$ and such
that the corresponding $x_{t}^{\pi}\geq0,$ for all $t\in\left[  0,T\right]  $.
We denote by $\Pi$ the class of such strategies.

Note that the cost functional $\left(  4.4\right)  $ may be rewritten in
men-field terms as%
\[
J\left(  \pi\right)  =\mathbb{E}\left[  -x_{T}+\mu\left(  x_{T}+\mathbb{E}%
\left[  -x_{T}\right]  \right)  ^{2}\right]  .
\]

\bigskip

2) \textbf{Example 2 : }Mean-field-type game with one risk-sensitive decision-maker

The mean-field state equation is given by

\begin{center}
$\left\{
\begin{array}
[c]{l}%
dX_{t}=\overline{b}(t,X_{t},P_{X_{t}},u_{t})dt+\sigma(t,X_{t},P_{X_{t}}%
)dW_{t}\\
X_{0}=x.
\end{array}
\right.  $
\end{center}

the drift term $\overline{b}$ has the special form $\overline{b}$ $=\left(
\int|b|^{\alpha}(.,t,X_{t},y,u(t))P_{X_{t}}(t,dy)\right)  ^{\frac{1}{\alpha}}$

$\bigskip$The risk-sensitive control problem is to minimize the following cost functional

$J^{\theta}(u($\textperiodcentered$))$ $=\dfrac{1}{\theta}\log E\left(
\exp\theta\left[
{\displaystyle\int\limits_{0}^{T}}
h(t,X_{t},P_{X_{t}},u_{t})dt+g(X_{T},P_{X_{T}})\right]  \right)  $

This kind of MFSDE models has been used to understand muscle contraction.
Other similar models have been widely studied in chemical kinetics,
statistical mechanics and economics, \cite{Temb} for further details.

\subsection{The relaxed control problem}

As it is well known in control theory that in the absence of convexity
conditions, an optimal control may fail to exist in the set $\mathcal{U}_{ad}$
of strict controls (see e.g. \cite{Fl, Maz, Val}). To be convinced let us
consider the following examples.

\textbf{Example 1}.

Minimize $J(u)=\int\nolimits_{0}^{1}\left(  X(t)\right)  ^{2}dt$ \ over the
set $\mathcal{U}_{ad}$ of measurable functions $u:[0,1]\rightarrow\{-1,1\}$,
where $X^{u}(t)$ is the solution of $dX(t)=u(t)dt,$ $X(0)=0.$ We have
$\inf_{u\in\mathcal{U}_{ad}}J(u)=0$.

Indeed, consider the sequence of Rademacher functions:

\begin{center}%
\[
u_{n}(t)=(-1)^{k}\text{ if }\frac{k}{n}\leq t\leq\frac{(k+1)}{n},0\leq k\leq
n-1.
\]

\end{center}

Then clearly $|X^{u_{n}}(t)|\leq1/n$ and $|J(u_{n})|\leq1/n^{2}$ which implies
that $\inf_{u\in\mathcal{U}_{ad}}J(u)=0$. There is however no control
$\widehat{u}$ such that $J(\widehat{u})=0$. If this were the case, then for
every $t$, $X^{\widehat{u}}(t)=0$. This in turn would imply that $\widehat
{u}_{t}=0$, which is impossible.

------------------------------------

\textbf{Example 2. }Minimize $J(u)=E\left[  \int\nolimits_{0}^{T}\left[
X^{u}(t)^{2}+(1-u(t))^{2}\right]  dt\right]  $subject to $dX^{u}%
(t)=u(t)dt+dW$, $X^{u}(0)=0,$ where $W$ is one dimensional Brownian motion,
$u$ is an open-loop control that is a measurable function from $\left[
0,T\right]  $ into $\left[  -1,1\right]  .$ The optimal control minimizes the
functional $\int\nolimits_{0}^{T}\left[  \widehat{X}^{u}(t)^{2}+(1-u(t))^{2}%
\right]  dt$ where $\widehat{X}^{u}(t)=E\left[  X^{u}(t)\right]  .$

It is not difficult to prove that the family of Rademacher functions is a
minimizing sequence ($\lim$ $J(u_{n})=0$), then it follows that $\inf
_{u\in\mathcal{U}_{ad}}J(u)=0$. But there is no control $u\in\mathcal{U}_{ad}$
satisfying $J(u)=0$, since it would have to satisfy $\widehat{X}^{u}(t)=0$ and
$\left\vert u(t)\right\vert $ $=1$ $a.e.$ at the same time, which is impossible.

The problem in both examples is that the sequence $(u_{n})$ has no limit in
the space of strict controls. This limit, if it exists, would be the natural
candidate for optimality.

These examples suggest that the set\ of strict controls is too narrow and
should be embedded into a wider class of relaxed controls, with nice
compactness properties. For the relaxed model, to be a true extension of the
original control problem, the following both conditions must be satisfied:

i) The value functions of the original and the relaxed control problems must
be equal.

ii) The relaxed control problem must have an optimal solution.

\bigskip

The idea of relaxed control is to replace the $\mathbb{A}$-valued process
$(u_{t})$ with a $\mathbb{P}(\mathbb{A})$-valued process $(\mu_{t})$, where
$\mathbb{P}(\mathbb{A})$ is the space of probability measures equipped with
the topology of weak convergence. Then $(\mu_{t})$ may be identified as a
random product measure on $[0,T]\times\mathbb{A}$, whose projection on $[0,T]$
coincides with Lebesgue measure.

Let $\mathbb{V}\ $be the set of product measures $\mu$ on $[0,T]\times
\mathbb{A}$ whose projection on $\left[  0,T\right]  $ coincides with the
Lebesgue measure $dt$. It is clear that every $\mu$ in $\mathbb{V}$ may be
disintegrated as $\mu=dt.\mu_{t}(da)$, where $\mu_{t}(da)$ is a transition
probability. The elements of $\mathbb{V}$ are called Young measures in
deterministic theory, see \cite{Val}.

$\mathbb{V}$ as a closed subspace of the space of positive Radon measures
$\mathbb{M}_{+}([0,T]\times\mathbb{A})$ is compact for the topology of weak
convergence. In fact it can be proved that it is compact also for the topology
of stable convergence, where test functions are measurable, bounded functions
$f(t,a)$ continous in $a,$ see \cite{EKNJ, JM} for further details.

\begin{definition}
A relaxed control on the filtered probability space $\left(  \Omega
,\mathcal{F},\mathcal{F}_{t},P\right)  $ is a random variable $\mu=dt.\mu
_{t}(da)$ with values in $\mathbb{V}$, such that $\mu_{t}(da)$ is
progressively measurable with respect to $(\mathcal{F}_{t})$ and such that for
each $t$, $1_{(0,t]}.\mu$ is $\mathcal{F}_{t}-$measurable. Let us denote
$\mathcal{R}$ the set of relaxed controls.
\end{definition}

\begin{remark}
The set $\mathcal{U}_{ad}$ of strict controls is embedded into the set
$\mathcal{R}$ of relaxed controls by identifying $u_{t}$ with $dt\delta
_{u_{t}}(da).$
\end{remark}

\bigskip

Let us come back to the first example. If we identify $u_{n}(t)$ with the
Dirac measure $\delta_{u_{n}(t)}(du),$ then it is not difficult to prove that
the sequence of product measures $(dt\delta_{u_{n}(t)}(du))_{n}$ converges
weakly to $(dt/2)$\textperiodcentered$\lbrack\delta_{-1}+\delta_{1}](da)$.

\bigskip

Let us define the relaxed model by%

\[
x_{t}^{q}=x_{0}+\int_{0}^{t}ds\int_{U}uq\left(  s,da\right)
\]
and the associated relaxed cost is given by%
\[
J\left(  q\right)  =\int_{0}^{1}\left(  x_{t}^{q}\right)  ^{2}dt
\]
Then it is clear that the strict control problem is generalized by the relaxed
problem, by taking measures $q$ of the form
\[
q\left(  dt,du\right)  =\delta_{u_{t}}\left(  du\right)  dt
\]
Moreover if
\[
\widehat{q}\left(  dt,du\right)  =\frac{1}{2}\left[  \delta_{-1}+\delta
_{1}\right]  \left(  du\right)  dt
\]
then we have$\ J\left(  \widehat{q}\right)  =0$ and $\widehat{q}$ as an
optimal relaxed control. Moreover since $\inf_{u\in\mathcal{U}_{ad}}%
J(u)=\inf_{q\in\mathcal{R}}J(q)=0,$ then the value functions of the strict and
relaxed control problems are the same.

\subsubsection{Discussion on the relaxation of the state process}

A natural question arises: what is the natural SDE associated with a relaxed
control. Let us point out that in the deterministic case or in the stochastic
case where only the drift is controlled, one has just to replace in equation
(\ref{MFSDE1}) the drift by the same drift integrated against the relaxed
control. Now we are in a situation where both the drift and diffusion
coefficient are controlled. Let us try a direct relaxation of the original
equation as in \cite{Cha}. \ The "relaxed" control problem will be governed by
the MFSDE

\begin{center}
$\left\{
\begin{array}
[c]{l}%
dX_{t}=%
{\displaystyle\int\limits_{A}}
b(t,X_{t},E(\Psi(X_{t})),a)\mu_{t}(da)dt+%
{\displaystyle\int\limits_{A}}
\sigma(t,X_{t},E(\Psi(X_{t})),a)\mu_{t}(da)dW_{t}\\
X_{0}=x
\end{array}
\right.  $
\end{center}

As it will be shown in the sequel, this model does not fullfill the
requirements of a true relaxed model. The reason is that the relaxed process
is not continuous in the control variable and as a consequence, the value
functions of the original and relaxed control problems are not equal. Let us
consider the following example.

Consider the control problem governed by the following SDE without mean-field terms:

\begin{center}
$\left\{
\begin{array}
[c]{l}%
dX_{t}=u_{t}dW_{t}\\
X_{0}=x
\end{array}
\right.  $
\end{center}

where the control $u\in\mathcal{U}_{ad}:$ the set of measurable functions
$u:[0,1]\rightarrow\mathbb{A=}\left[  -1,1\right]  .$

The "relaxed" model will be governed by the equation

\begin{center}
$\left\{
\begin{array}
[c]{l}%
dX_{t}=%
{\displaystyle\int\nolimits_{\mathbb{A}}}
a\mu_{t}(da)dW_{t}\\
X_{0}=x
\end{array}
\right.  $
\end{center}

Consider the sequence of Rademacher functions

\begin{center}
$u_{n}(t)=(-1)^{k}$ if $\frac{k}{n}\leq t\leq\frac{(k+1)}{n}$, $0\leq k\leq
n-1.$
\end{center}

\begin{lemma}
Let $dt.\delta_{u_{n}(t)}(da)$ the relaxed control associated to $u_{n}(t),$
then \ the sequence $\left(  dt.\delta_{u_{n}(t)}(da)\right)  $ converges
weakly to $dt(1/2)(\delta_{-1}+\delta_{1})(da).$
\end{lemma}

\bop See \cite{Maz} Lemma 1.1, page 20 \eop

\begin{remark}
The sequence of Rademacher functions is a typical example of a minimising
sequence with no limit in the set of strict controls. However its weak limit
$dt(1/2)(\delta_{-1}+\delta_{1})(da)$ is an optimal relaxed control \cite{Maz,
Val}.
\end{remark}

It is clear that $X_{t}^{n}=%
{\displaystyle\int\limits_{0}^{t}}
u_{n}(s)$.$dW_{s}=%
{\displaystyle\int\limits_{0}^{t}}
\left[
{\displaystyle\int\nolimits_{\mathbb{A}}}
a\delta_{u_{n}(s)}(da)\right]  dW_{s}$ is a continuous martingale with
quadratic variation $\left\langle X^{n},X^{n}\right\rangle _{t}=%
{\displaystyle\int\limits_{0}^{t}}
u_{n}^{2}(s)$.$ds=t.$ Therefore $\left(  X_{t}^{n}\right)  $ is a Brownian motion.

By Lemma 2.4, the sequence of relaxed controls $\left(  dt.\delta_{u_{n}%
(t)}(da)\right)  $ converges weakly in $\mathbb{M}_{+}([0,T]\times\mathbb{A})$
to $\mu^{\ast}=(1/2)dt(\delta_{-1}+\delta_{1})(da).$ Let $X^{\ast}$ be the
relaxed state process corresponding to the limit $\mu^{\ast},$ then

\begin{center}
$X^{\ast}(t)=%
{\displaystyle\int\limits_{0}^{t}}
{\displaystyle\int\limits_{\mathbb{A}}}
a.(1/2)(\delta_{-1}+\delta_{1})(da)dW_{t}=0.$
\end{center}

It is obvious that the sequence of state processes $\left(  X_{t}^{n}\right)
$ cannot converge in $L^{2}$ to $X_{t}^{\ast}.$ Indeed

\begin{center}
$E\left[  \left\vert X_{t}^{n}-X_{t}^{\ast}\right\vert ^{2}\right]  =E\left[
\left\vert X_{t}^{n}\right\vert ^{2}\right]  =E\left[  \left\vert
{\displaystyle\int\limits_{0}^{t}}
u_{n}(s).dW_{s}\right\vert ^{2}\right]  =%
{\displaystyle\int\limits_{0}^{t}}
u_{n}^{2}(s)$.$ds=t$
\end{center}

This implies in particular that the state process is not continuous in the
control variable and as a byproduct, the value functions of the strict and
"relaxed" control problems are not equal. Moreover, even if the set
$\mathbb{V}$ is compact, there is no mean to prove existence of an optimal
control for this model.

\textbf{What is the right relaxed state process?}

The reason why the proposed model in \cite{Cha} is not a true extension of the
original strict control problem, is that the stochastic integral part does not
behave as a Lebesgue integral. In fact, one should relax the drift and the
quadratic variation of the martingale part, which is a Lebesgue integral.

In the relaxed model, the quadratic variation process must be $%
{\displaystyle\int\limits_{0}^{t}}
{\displaystyle\int\limits_{A}}
\sigma\sigma^{\ast}(t,X_{t},E(\Phi(X_{t})),a)\mu_{t}(da)dt,$ which is more
natural than relaxing the stochastic integral itself.$\ $Now, one has to look
for a square integrable martingale whose quadratic variation is given by $%
{\displaystyle\int\limits_{0}^{t}}
{\displaystyle\int\limits_{\mathbb{A}}}
\sigma\sigma^{\ast}(t,X_{t},E(\Phi(X_{t}),a)\mu_{t}(da)dt,$ which is
equivalent to the search of an object which is a martingale whose quadratic
variation is $dt\mu_{t}(da).$ This object is precisely a continuous orthogonal
martingale measure, whose covariance measure is $dt\mu_{t}(da).$ This is
equivalent to the relaxation of the infinitesimal generator associated to the
state process.

Following \cite{JMW} Prop. 1.10, existence of a weak solution of equation
(\ref{MFSDE1}) is equivalent to existence of a solution for the non linear
martingale problem:

\begin{center}
$f(X_{t})-f(X_{0})-\int\limits_{0}^{t}L^{P_{X_{s}}}f(s,X_{s},u_{s})\,ds$
\textit{is a} $P-$\textit{martingale}$\mathit{,}$
\end{center}

\textit{for every }$f\in C_{b}^{2},$\textit{\ for each }$t>0,$\textit{\ where
}$L$\textit{\ is the infinitesimal generator }associated with equation
(\ref{MFSDE1}),
\[
L^{\nu}f(t,x,a)=\frac{1}{2}\sum\limits_{i,j}\left(  a_{i,j}\frac{\partial
^{2}f}{\partial x_{i}\partial x_{j}}\right)  (t,x,a)+\sum\limits_{i}\left(
b_{i}\frac{\partial f}{\partial x_{i}}\right)  (t,x,a),
\]

$b=b(t,x,\left\langle \Psi,\nu\right\rangle ,a)$ and $a_{i,j}=\sigma
\sigma^{\ast}(t,x,\left\langle \Phi,\nu\right\rangle ,a),$ $\nu\in
\mathbb{M}_{1}(\mathbb{R}^{d}).$

The natural relaxed martingale problem associated to the relaxed generator is
defined as follows:

\begin{center}%
\begin{equation}
f(X_{t})-f(X_{0})-\int\limits_{0}^{t}%
{\displaystyle\int\limits_{\mathbb{A}}}
L^{P_{X_{s}}}f(s,X_{s},a)\,\mu_{s}(da)ds\text{ }\mathit{is\ a}\text{
}P-\mathit{martingale} \label{NLRM}%
\end{equation}

\end{center}

\textit{for each }$f\in C_{b}^{2},$\textit{\ for each }$t>0.$

Now, what is the stochastic differential equation corresponding to the relaxed
martingale problem (\ref{NLRM})? The answer is given by the following theorem.

\begin{theorem}
\ \textit{1)Let }$P$\textit{\ be the solution of the martingale problem
(\ref{NLRM}). Then }$P$\textit{\ is the law of a }$d$\textit{--dimensional
adapted and continuous process }$X$\textit{\ defined on an extension of the
space }$\left(  \Omega,\mathcal{F},\mathcal{F}_{t},P\right)  $\textit{\ and
solution of the following MFSDE starting at }$x$\textit{:}
\end{theorem}

\begin{equation}
\left\{
\begin{array}
[c]{l}%
dX_{t}=\int\nolimits_{\mathbb{A}}b(t,X_{t},E\left(  \Psi\left(  X_{t}\right)
\right)  a)\,\mu_{t}(da)dt+\int\nolimits_{\mathbb{A}}\sigma(t,X_{t}%
,E(\Phi(X_{t}),a)\,M(da,dt),\\
X_{0}=x
\end{array}
\right.  \label{RMFSDE}%
\end{equation}

\textit{where }$M=(M^{k})_{k=1}^{d}$\textit{\ is a family of }$d$%
\textit{-strongly orthogonal continuous martingale measures, each having
intensity }$\mathit{dt}\mu_{t}(da).$

\textit{2) If the coefficients }$b$\textit{\ and }$\sigma$\textit{\ are
Lipschitz in }$x$\textit{, }$\mathit{y}$\textit{, uniformly in }%
$t$\textit{\ and }$a$\textit{, the SDE (2.6}$)$\textit{\ has a unique pathwise
solution.}

\bop1) The proof is based essentially on the same arguments as in \cite{EM}
Theorem IV-2 and \cite{JMW} Prop. 1.10.

2) The coefficients being Lipschitz continuous, following the same steps as in
\cite{JMW} and \cite{EM}, it is not difficult to prove that Equation
(\ref{RMFSDE}) has a unique solution such that for every $p>0$ we have
$E(\left\vert X_{t}\right\vert ^{p})<+\infty.$ \eop

\begin{remark}
Note that the othogonal martingale measure corresponding to the relaxed
control $\mathit{dt}\mu_{t}(da)$ is not unique.
\end{remark}

Let us give the precise definition of a martingale measure introduced by Walsh
\cite{Wal}, see also \cite{EM, Mel} for more details.

\begin{definition}
Let $(\Omega,\mathcal{F},\mathcal{F}_{t},P)$ be a filtered probability space,
and $M(t,B)$ a random process, where $B\in\mathcal{B}\left(  \mathbb{A}%
\right)  .$ $M$ is a ($\mathcal{F}_{t},P)-$martingale measure if:

1)For each $B\in\mathcal{B}\left(  \mathbb{A}\right)  ,\left(  M(t,B)\right)
_{t\geq0}$ is a square integrable martingale, $M(0,B)=0$.

2)For every $t>0$, $M(t,.)$ is a $\sigma-$finite $L^{2}$-valued measure.

It is called continuous if for each $B\in\mathcal{B}\left(  \mathbb{A}\right)
,$ $M(t,B)$ is continuous and orthogonal if $M(t,B).M(t,C)$ is a martingale
whenever $B\cap C=\phi.$
\end{definition}

\begin{remark}
When the martingale measure $M$ is orthogonal, it is proved in \cite{Wal} the
existence of a random positive $\sigma$-finite measure $\mu\left(
dt,da\right)  \ $on $\left[  0,T\right]  \times\mathbb{A}$ such that
$\left\langle M(.,B),M(.,B)\right\rangle _{t}=\mu\left(  \left[  0,t\right]
\times B\right)  $ for all $t>0$ and $B\in\mathcal{B}\left(  \mathbb{A}%
\right)  .$ $\mu\left(  dt,da\right)  $ is called the covariance measure of
$M$.
\end{remark}

\textbf{Example \thinspace}Let ${\mathbb{A}}=\left\{  a_{1},a_{2},\cdots
,a_{n}\right\}  $ be a finite set. Then every relaxed control $dt\,\mu
_{t}(da)$ will be a convex combination of the Dirac measures $dt\,\mu
_{t}(da)=\sum_{i=1}^{n}\alpha_{t}^{i}\,dt\,\delta_{a_{i}}(da),$ where for each
$i$, $\alpha_{t}^{i}$ is a real--valued progressively measurable process, such
that $0\leq\alpha_{t}^{i}\,\leq1$ and $\sum_{i=1}^{n}\alpha_{t}^{i}\,=1$. It
is obvious that the solution of the relaxed martingale problem
\textit{\ref{NLRM}} is the law of the solution of the SDE
\begin{equation}
dX_{t}=\sum\limits_{i=1}^{n}b(t,X_{t},E\left(  \Psi\left(  X_{t}\right)
\right)  ,a_{i})\alpha_{t}^{i}dt+\sum\limits_{i=1}^{n}\sigma(t,X_{t},E\left(
\Psi\left(  X_{t}\right)  \right)  ,a_{i})\left(  \alpha_{t}^{i}\right)
^{1/2}\,dB_{t}^{i},\quad\quad X_{0}=x, \label{FinSDE}%
\end{equation}
where the $B^{i}$'s are independent $d$-dimensional Brownian motions, on an
extension of the initial probability space. The process $M$ defined by
\[
M(A\times\left[  0,t\right]  )=\sum\limits_{i=1}^{n}\int\limits_{0}^{t}\left(
\alpha_{s}^{i}\right)  ^{1/2}1_{\left\{  a_{i}\in A\right\}  }dB_{s}^{i}%
\]
is in fact an orthogonal continuous martingale measure (cf. \cite{EKNJ, Wal})
with intensity $\mu_{t}(da)dt=\sum\alpha_{t}^{i}\,\delta_{a_{i}}(da)dt$. Thus,
the SDE (\ref{FinSDE}) can be expressed in terms of $M$ and $\mu$ as follows:%

\[
dX_{t}=\int\nolimits_{\mathbb{A}}b(t,X_{t},E\left(  \Psi\left(  X_{t}\right)
\right)  a)\,\mu_{t}(da)dt+\int\nolimits_{\mathbb{A}}\sigma(t,X_{t}%
,E(\Phi(X_{t})),a)\,M(da,dt)
\]

\subsubsection{Approximation of the relaxed model}

The relaxed control problem is now driven by equation%

\begin{equation}
\left\{
\begin{array}
[c]{l}%
dX_{t}=\int\nolimits_{\mathbb{A}}b(t,X_{t},E\left(  \Psi\left(  X_{t}\right)
\right)  a)\,\mu_{t}(da)dt+\int\nolimits_{\mathbb{A}}\sigma(t,X_{t}%
,E(\Phi(X_{t})),a)\,M(da,dt),\\
X_{0}=x
\end{array}
\right.  \label{RMFSDE2}%
\end{equation}

and accordingly the relaxed cost functional is given by%

\begin{equation}
J(\mu)=E\left(
{\displaystyle\int\limits_{0}^{T}}
\int\nolimits_{\mathbb{A}}h(t,X_{t},E(\varphi(X_{t})),a)\mu_{t}(da)dt+g(X_{T}%
,E(\lambda(X_{T}))\right)  .
\end{equation}

We show in this section that the strict and the relaxed control problems have
the same value function. This is based on the chattering lemma and the
stability of the state process with respect to the control variable.

\begin{lemma}
\textbf{(}Chattering lemma\textbf{)\ }i)\textbf{ }\textit{Let }$(\mu_{t}%
)$\textit{\ be a relaxed control}$.$\textit{\ Then there exists a sequence of
adapted processes }$(u_{t}^{n})$\textit{\ with values in }$\mathbb{A}%
$\textit{, such that the sequence of random measures }$\left(  \delta
_{u_{t}^{n}}(da)\,dt\right)  $\textit{\ converges in }$\mathbb{V}$\textit{ to
}$\mu_{t}(da)\,dt,$\textit{\ }$P-a.s.$

ii) For any $g$ continuous in $\left[  0,T\right]  \times\mathbb{M}%
_{1}(\mathbb{A})$ such that $g(t,.)$ is linear, we have $P-a.s$%
\begin{equation}
\underset{n\rightarrow+\infty}{\lim}%
{\displaystyle\int\limits_{0}^{t}}
g(s,\delta_{u_{s}^{n}})ds=%
{\displaystyle\int\limits_{0}^{t}}
g(s,\mu_{s})ds\text{ uniformly in }t\in\left[  0,T\right]  .
\label{Uniform Limit}%
\end{equation}

\end{lemma}

\bop See \cite{EKNJ} and \cite{Fl} Lemma 1 page 152.\eop

\begin{proposition}
1) Let $\mu=\mu_{t}(da)\,dt$ a relaxed control. Then there exist a continuous
orthogonal martingale measure $M(dt,da)$ whose covariance measure is given by
$\mu_{t}(da)\,dt.$

2) If we denote $M^{n}(t,B)=\int\nolimits_{0}^{t}\int\nolimits_{B}%
\delta_{u_{s}^{n}}(da)dW_{s},$ where $\left(  u^{n}\right)  $ is defined as in
the last Lemma, then for every bounded predictable process $\varphi
:\Omega\times\left[  0,T\right]  \times\mathbb{A}\rightarrow\mathbb{R}$, such
that $\varphi(\omega,t,.)$ is continuous$,$ we have

$E\left[  \left(  \int\nolimits_{0}^{t}\int\nolimits_{\mathbb{A}}%
\varphi(\omega,t,a)M^{n}(dt,da)-\int\nolimits_{0}^{t}\int\nolimits_{\mathbb{A}%
}\varphi(\omega,t,a)M(dt,da)\right)  ^{2}\right]  \rightarrow0$ as
$n\longrightarrow+\infty,$

for a suitable Brownian motion $W$ defined on an eventual extension of the
probability space.
\end{proposition}

\begin{center}

\end{center}

\bop See \cite{Mel} pages 196-197.\eop

\begin{proposition}
\textit{1) Let }$X_{t},$ $X_{t}^{n}$ \textit{be the solutions of state
equation (}\ref{RMFSDE}) corresponding to $\mu$ and $u^{n},$ where $\mu$ and
$u^{n}$ are defined as in the chattering lemma\textit{. Then \ }%
\begin{equation}
\underset{n\rightarrow\infty}{\lim}E\left[  \underset{0\leq t\leq T}{\sup
}\left\vert X_{t}^{n}-X_{t}\right\vert ^{2}\right]  =0.
\end{equation}

2) Let $J(u^{n})$ and $J(\mu)$ be the expected costs corresponding
respectively to $u^{n}$ and $\mu,$ then $\left(  J\left(  u^{n}\right)
\right)  $ converges to $J\left(  \mu\right)  .$
\end{proposition}

\bop

1)Let $\mu$ a relaxed control and $\left(  dt\delta_{u_{t}^{n}}(da)\right)  $
the sequence of atomic measures associated to the sequence of strict controls
$\left(  u^{n}\right)  ,$ as in the last Lemma. Let $X_{t},$ $X_{t}^{n}$ the
corresponding state processes. If we denote $M^{n}(t,B)=\int\nolimits_{0}%
^{t}\int\nolimits_{B}\delta_{u_{s}^{n}}(da)dW_{s},$ then $X^{n}$ may be
written in a relaxed form

\begin{center}
$\left\{
\begin{array}
[c]{l}%
dX_{t}^{n}=%
{\displaystyle\int\limits_{\mathbb{A}}}
b(t,X_{t}^{n},E(\Psi(X_{t}^{n})),a)\delta_{u_{t}^{n}}(da)dt+%
{\displaystyle\int\limits_{\mathbb{A}}}
\sigma(t,X_{t},E(\Phi(X_{t})),a)M^{n}(dt,da)\\
X_{0}=x
\end{array}
\right.  $
\end{center}

We have%

\[%
\begin{array}
[c]{cl}%
\left\vert X_{t}-X_{t}^{n}\right\vert  & \leq\left\vert \int\nolimits_{0}%
^{t}\int\nolimits_{\mathbb{A}}b\left(  s,X_{s},E(\Psi(X_{s})),u\right)
\mu_{s}(da).ds-\int\nolimits_{0}^{t}\int\nolimits_{\mathbb{A}}b\left(
s,X_{s}^{n},E(\Psi(X_{s}^{n})),u\right)  \delta_{u_{s}^{n}}(da)ds\right\vert
\\
& +\left\vert \int\nolimits_{0}^{t}\int\nolimits_{\mathbb{A}}\sigma\left(
s,X_{s},E(\Phi(X_{s})),a\right)  M(ds,da)-\int\nolimits_{0}^{t}\int
\nolimits_{\mathbb{A}}\sigma\left(  s,X_{s}^{n},E(\Phi(X_{s}^{n})),a\right)
M^{n}(ds,da)\right\vert \\
& \leq\left\vert \int\nolimits_{0}^{t}\int\nolimits_{\mathbb{A}}b\left(
s,X_{s},E(\Psi(X_{s})),u\right)  \mu_{s}(da).ds-\int\nolimits_{0}^{t}%
\int\nolimits_{\mathbb{A}}b\left(  s,X_{s},E(\Psi(X_{s})),a\right)
\delta_{u_{s}^{n}}(da)ds\right\vert \\
& +\left\vert \int\nolimits_{0}^{t}\int\nolimits_{\mathbb{A}}b\left(
s,X_{s},E(\Psi(X_{s})),u\right)  \delta_{u_{s}^{n}}(da).ds-\int\nolimits_{0}%
^{t}\int\nolimits_{\mathbb{A}}b\left(  s,X_{s}^{n},E(\Psi(X_{s}^{n}%
)),a\right)  \delta_{u_{s}^{n}}(da)ds\right\vert \\
& +\underset{s\leq t}{\sup}\left\vert \int\nolimits_{0}^{s}\int
\nolimits_{\mathbb{A}}\sigma\left(  v,X_{v},E(\Phi(X_{v})),a\right)
M(dv,da)-\int\nolimits_{0}^{t}\int\nolimits_{\mathbb{A}}\sigma\left(
v,X_{v},E(\Phi(X_{v})),a\right)  M^{n}(dv,da)\right\vert \\
& +\underset{s\leq t}{\sup}\left\vert \int\nolimits_{0}^{s}\int
\nolimits_{\mathbb{A}}\sigma\left(  v,X_{v},E(\Phi(X_{v}),a\right)
M^{n}(dv,da)-\int\nolimits_{0}^{t}\int\nolimits_{\mathbb{A}}\sigma\left(
v,X_{v}^{n},E(\Phi(X_{v}^{n}),a\right)  M^{n}(dv,da)\right\vert
\end{array}
\]

Then by using Burkholder-Davis-Gundy inequality for the martingale part and
the fact that all the functions in equation \textit{(}\ref{RMFSDE}) are
Lipschitz continuous, it holds that

\begin{center}%
\begin{equation}
E\left(  \underset{0\leq t\leq T}{\sup}\left\vert X_{t}-X_{t}^{n}\right\vert
^{2}\right)  \leq K\left[  \int\nolimits_{0}^{T}E\left(  \underset{0\leq s\leq
t}{\sup}\left\vert X_{s}-X_{s}^{n}\right\vert ^{2}\right)  dt+\varepsilon
_{n}\right]  ,
\end{equation}

\end{center}

where $K$ is a nonnegative constant and

\begin{center}%
\begin{multline}
\varepsilon_{n}=E\left(  \underset{0\leq t\leq T}{\sup}\left\vert
\int\nolimits_{0}^{t}\int\nolimits_{\mathbb{A}}b\left(  s,X_{s},E(\Psi
(X_{s})),u\right)  \mu_{s}(da)ds-\int\nolimits_{0}^{t}\int
\nolimits_{\mathbb{A}}b\left(  s,X_{s},E(\Psi(X_{s})),a\right)  \delta
_{u_{s}^{n}}(da)ds\right\vert ^{2}\right) \\
+E\left(  \underset{0\leq t\leq T}{\sup}\left\vert \int\nolimits_{0}^{t}%
\int\nolimits_{\mathbb{A}}\sigma\left(  s,X_{s},E(\Psi(X_{s})),a\right)
M(ds,da)-\int\nolimits_{0}^{t}\int\nolimits_{\mathbb{A}}\sigma\left(
s,X_{s},E(\Psi(X_{s})),a\right)  M^{n}(ds,da)\right\vert ^{2}\right)  .
\end{multline}

\end{center}

By using Lemma 2.10 and the Lebesgue dominated convergence theorem, it holds
that $\underset{n\rightarrow+\infty}{\lim}\varepsilon_{n}=0.$ We conclude by
using Gronwall's Lemma.

2) Property$\ $1) implies that the sequence $\left(  X_{t}^{n}\right)  $
converges to $X_{t}$ in probability uniformly in $t,$ then we have%

\[%
\begin{array}
[c]{cl}%
\left\vert J\left(  u^{n}\right)  -J\left(  \mu\right)  \right\vert  & \leq
E\left[  \int\limits_{0}^{T}\int_{\mathbb{A}}\left\vert h(t,X_{t}%
^{n},E(\varphi(X_{t}^{n})),a)-\,h(t,X_{t},E(\varphi(X_{t})),a)\right\vert
\delta_{u_{t}^{n}}(da)\,dt\right] \\
& +E\left[  \left\vert \int\limits_{0}^{T}\int_{\mathbb{A}}h(t,X_{t}%
,E(\varphi(X_{t})),a)\delta_{u_{t}^{n}}(da)\,dt-\int\limits_{0}^{T}%
\int_{\mathbb{A}}h(t,X_{t},E(\varphi(X_{t})),a)\mu_{t}(da)\,dt\right\vert
\right] \\
& +E\left[  \left\vert g(X_{T}^{n},E(\lambda(X_{T}^{n}))-g(X_{T}%
,E(\lambda(X_{T}))\right\vert \right]
\end{array}
\]

It follows from the continuity and boundness of the functions $h$, $g$,
$\varphi$ and $\lambda$ with respect to $x$ and $y,$ that the first and third
terms in the right hand side converge to $0$ . The second term in the right
hand side tends to $0$ by the weak convergence of the sequence $\mu^{n}$ to
$\mu,$ the continuity and the boundedness of $h$ in the variable $a$. We use
the dominated convergence theorem to conclude. \eop

\begin{remark}
As a consequence of the last proposition, it holds that the infimum among
relaxed controls is equal to the infimum among strict controls, that is the
relaxed model is a true extension of the original control problem.
\end{remark}

\subsection{Existence of optimal relaxed controls}

The following theorem which is the main result of this section, extends
\cite{BahDjeMez, EKNJ, Fl} to systems driven by mean field SDEs with
controlled diffusion coefficient.

\begin{theorem}
Under assumptions $\mathbf{(H_{1})}$, $\mathbf{(H_{2})}$, there exist an
optimal relaxed control.
\end{theorem}

The proof is based on some auxiliary Lemmas and will be given later.

Let $\left(  \mu^{n}\right)  _{n\geq0}$ be a minimizing sequence, that is
$\underset{n\rightarrow\infty}{\lim}J\left(  \mu^{n}\right)  =\underset
{q\in\mathcal{R}}{\inf}J\left(  \mu\right)  $ and let $X^{n}$ be the unique
solution of (\ref{RMFSDE}), associated with $\mu^{n}$ and $M^{n}$ where
$M^{n}$ is a continuous orthogonal martingale measure with intensity $\mu
^{n}.$ We will prove that the sequence $(\mu^{n},M^{n},X^{n})$ is tight and
then show that we can extract a subsequence, which converges in law to a
process $(\widehat{\mu},\widehat{M},\widehat{X}),$ which satisfies the same
MFSDE. To finish the proof we show that the sequence of cost functionals
$(J(\mu^{n}))_{n}$ converges to $J(\widehat{\mu})$ which is equal to
$\underset{\mu\in\mathcal{R}}{\inf}J\left(  \mu\right)  $ and then we conclude
that $(\widehat{\mu},\widehat{M},\widehat{X})$ is optimal.

\begin{lemma}
The sequence of distributions of the relaxed controls $(\mu^{n})_{n}$ is
relatively compact in $\mathbb{V}$.
\end{lemma}

\bop The relaxed controls $\mu^{n}$ are random variables with values in the
space $\mathbb{V}$ which is compact. Then Prohorov's theorem yields that the
family of distributions associated to $(\mu^{n})_{n\geq0}$ is tight then it is
relatively compact. \eop

\begin{lemma}
\textit{The family of martingale measures }$\left(  M^{n}\right)  _{n\geq0}%
$\textit{\ is tight in the space }$C\left(  \left[  0,1\right]  ,\mathcal{S}%
^{\prime}\right)  $\textit{\ of continuous functions from }$\left[
0,1\right]  $\textit{\ into }$\mathcal{S}^{\prime},$\textit{\ the topological
dual of the Schwartz space }$\mathcal{S}$\textit{\ of rapidly decreasing
functions}$.$
\end{lemma}

\bop The martingale measures $M^{n},$ $n\geq0,$ can be considered as random
variables with values in $C\left(  \left[  0,1\right]  ,\mathcal{S}^{\prime
}\right)  $. By \cite{Mit}, Therem 5.1, it is sufficient to show that for
every $\varphi$ in $\mathcal{S}$ the family $\left(  M^{n}\left(
\varphi\right)  ,n\geq0\right)  $ is tight in $C\left(  \left[  0,T\right]
,\mathbb{R}^{d}\right)  $ where $M^{n}\left(  \omega,t,\varphi\right)
=\int_{\mathbb{A}}\varphi(a)M^{n}(\omega,t,da).$ Let $p>1$ and $s<t$, by the
Burkholder-Davis-Gundy inequality we have
\begin{align*}
E\left(  \left\vert M_{t}^{n}(\varphi)-M_{s}^{n}(\varphi)\right\vert
^{2p}\right)   &  \leq C_{p}\,E\left[  \left(  \int\limits_{s}^{t}%
\int_{\mathbb{A}}\left\vert \varphi(a)\right\vert ^{2}\mu_{t}^{n}%
(da)\,dt\right)  ^{p}\right] \\
&  \leq C_{p}\underset{a\in\mathbb{A}}{\sup}\left\vert \varphi(a)\right\vert
^{2p}\left\vert t-s\right\vert ^{p}=K_{p}\,\left\vert t-s\right\vert ^{p},
\end{align*}
where $K_{p}$ is a constant depending on $p$ and $\varphi$. Then, the
Kolmogorov tightness criteria in $C\left(  \left[  0,T\right]  ,\mathbb{R}%
^{d}\right)  $ is fulfilled and the sequence $\left(  M^{n}\left(
\varphi\right)  \right)  $ is tight. \eop

\begin{lemma}
\textit{The sequence }$\left(  X^{n}\right)  _{n\geq0}$\textit{\ is tight in
the space }$C\left(  \left[  0,T\right]  ,\mathbb{R}^{d}\right)  $
\end{lemma}

\bop Let $p>1$ and $s<t$. Using usual arguments from stochastic calculus and
the boundness of the coefficients $b$ and $\sigma,$ it is easy to show that
\[
E\left(  \left\vert X_{t}^{n}-X_{s}^{n}\right\vert ^{2p}\right)  \leq
C_{p}\,\left\vert t-s\right\vert ^{p}%
\]
which yields the tightness of $\left(  X_{t}^{n},n\geq0\right)  $ \textit{in
}$C\left(  \left[  0,T\right]  ,\mathbb{R}^{d}\right)  $\eop

\bop\textbf{of Theorem 2.14}

By using the Lemmas 2.15, 2.16 and 2.17, it holds that the sequence of
processes $(\mu^{n},M^{n},X^{n})$ is tight on the space $\Gamma
=\mathbb{V\times}C\left(  \left[  0,1\right]  ,\mathcal{S}^{\prime}\right)
\times C\left(  \left[  0,T\right]  ,\mathbb{R}^{d}\right)  $. Then by the
Skorokhod representation theorem, there exists a probability space $\left(
\widehat{\Omega},\widehat{\mathcal{F}},\widehat{P}\right)  $, a sequence
$\widehat{\gamma}^{n}=\left(  \widehat{\mu}^{n},\widehat{M^{n}},\widehat
{X^{n}}\right)  $ and $\widehat{\gamma}=\left(  \widehat{\mu},\widehat
{M},\widehat{X}\right)  $ defined on this space such that:

(i) for each $n\in\mathbb{N}$, law$\left(  \gamma^{n}\right)  =$ law$\left(
\widehat{\gamma}^{n}\right)  $,

(ii) there exists a subsequence $\left(  \widehat{\gamma}^{n_{k}}\right)
$\ of $\left(  \widehat{\gamma}^{n}\right)  $, still denoted by $\left(
\widehat{\gamma}^{n}\right)  $, which converges to $\widehat{\gamma}%
,\widehat{P}$-$a.s$. on the space $\Gamma.$

This means in particular that the sequence of relaxed controls $(\widehat{\mu
}^{n})$ converges in the weak topology to $\widehat{\mu},$ $\widehat{P}-a.s.$
and $\left(  \widehat{M^{n}},\widehat{X^{n}}\right)  $ converges to $\left(
\widehat{M},\widehat{X}\right)  ,$ $\widehat{\mathbb{P}}-a.s.$ in $C\left(
\left[  0,1\right]  ,\mathcal{S}^{\prime}\right)  \times C\left(  \left[
0,T\right]  ,\mathbb{R}^{d}\right)  $.

According to property (i), we get

\begin{center}%
\begin{equation}
_{\left\{
\begin{array}
[c]{l}%
\widehat{X_{t}^{n}}=x+\int\nolimits_{0}^{t}\int\nolimits_{\mathbb{A}}b\left(
s,\widehat{X_{s}^{n}},E(\Psi(\widehat{X_{s}^{n}})),a\right)  \widehat{\mu}%
_{s}^{n}(da)ds+\int\nolimits_{0}^{t}\int\nolimits_{\mathbb{A}}\sigma\left(
s,\widehat{X_{s}^{n}},E(\Phi(\widehat{X_{s}^{n}})),a\right)  \widehat{M^{n}%
}(ds,da),\\
\widehat{X_{0}^{n}}=x.
\end{array}
\right.  }%
\end{equation}

\end{center}

Since the coefficients $b,$ $\sigma,$ $\Psi$ and $\Phi$ are Lipschitz
continuous in $(x,y),$then according to property (ii) and using similar
arguments as in \cite{Sk} page 32, it holds that%

\[
\int\nolimits_{0}^{t}\int\nolimits_{\mathbb{A}}b\left(  s,\widehat{X_{s}^{n}%
},E(\Psi(\widehat{X_{s}^{n}})),a\right)  \widehat{\mu}_{s}^{n}(da)ds\text{
converges in probability to }\int\nolimits_{0}^{t}\int\nolimits_{\mathbb{A}%
}b\left(  s,\widehat{X_{s}},E(\Psi(\widehat{X_{s}})),a\right)  \widehat{\mu
}_{s}(da)ds
\]

and%

\[
\int\nolimits_{0}^{t}\int\nolimits_{\mathbb{A}}\sigma\left(  s,\widehat
{X_{s}^{n}},E(\Phi(\widehat{X_{s}^{n}})),a\right)  \widehat{M^{n}%
}(ds,da)\text{ converges in probability to }\int\nolimits_{0}^{t}%
\int\nolimits_{\mathbb{A}}\sigma\left(  s,\widehat{X_{s}},E(\Phi
(\widehat{X_{s}})),a\right)  \widehat{M}(ds,da).
\]

$b$ and $\sigma$ are Lipschitz continuous, therefore $\widehat{X}$ is the
unique solution of the MFSDE

\begin{center}%
\begin{equation}
\left\{
\begin{array}
[c]{l}%
\widehat{X_{t}}=x+\int\nolimits_{0}^{t}\int\nolimits_{\mathbb{A}}b\left(
s,\widehat{X_{s}},E(\Psi(\widehat{X_{s}})),a\right)  \widehat{\mu}%
_{s}(da)ds+\int\nolimits_{0}^{t}\int\nolimits_{\mathbb{A}}\sigma\left(
s,\widehat{X_{s}},E(\Phi(\widehat{X_{s}})),a\right)  \widehat{M}(ds,da),\\
\widehat{X_{0}}=x.
\end{array}
\right.
\end{equation}

\end{center}

To finish the proof of Theorem 2.14, it remains to check that $\widehat{\mu}$
is an optimal control.

The functions $b$ and $\sigma$ are Lipschitz continuous, then according to the
above properties (i)-(ii) we get
\begin{align*}
\underset{\mu\in\mathcal{R}}{\inf}J\left(  \mu\right)   &  =\underset
{n\rightarrow\infty}{\lim}J\left(  \mu^{n}\right)  ,\\
&  =\underset{n\rightarrow\infty}{\lim}E\left[
{\displaystyle\int\limits_{0}^{T}}
\int\nolimits_{\mathbb{A}}h(t,X_{t}^{n},E(\varphi(X_{t}^{n})),a)\mu_{t}%
^{n}(da)dt+g(X_{T}^{n},E\lambda(X_{T}^{n}))\right] \\
&  =\underset{n\rightarrow\infty}{\lim}\widehat{E}\left[
{\displaystyle\int\limits_{0}^{T}}
\int\nolimits_{\mathbb{A}}h\left(  t,\widehat{X_{t}^{n}},E(\varphi
(\widehat{X_{t}^{n}})),a\right)  \widehat{\mu}_{t}^{n}(da)dt+g(\widehat
{X_{T}^{n}},E\lambda(\widehat{X_{T}^{n}}))\right] \\
&  =\widehat{E}\left[
{\displaystyle\int\limits_{0}^{T}}
\int\nolimits_{\mathbb{A}}h\left(  t,\widehat{X_{t}},E(\varphi(\widehat{X_{t}%
})),a\right)  \widehat{\mu}_{t}(da)dt+g(\widehat{X_{T}},E\lambda
(\widehat{X_{T}}))\right]  .
\end{align*}

Hence $\widehat{\mu}$ is an optimal control. \eop

\begin{remark}
The proof of the last Theorem is based on tightness and weak convergence
techniques. Then it is possible to prove it by using the non linear martingale
problem and following the same steps as in \cite{EKNJ} without using the
pathwise representation of the solution.
\end{remark}

We prove that in the next proposition that we can restrict the investigation
for an optimal relaxed control to the class of so-called sliding controls also
known as chattering controls, having the form
\begin{equation}
\mu_{t}=%
{\displaystyle\sum\limits_{i=1}^{p}}
\alpha_{i}(t)\delta_{u_{i}(t)}(da),u_{i}(t)\in\mathbb{A},\alpha_{i}%
(t)\geq0\text{ and }%
{\displaystyle\sum\limits_{i=1}^{p}}
\alpha_{i}(t)=1.
\end{equation}

where $\alpha_{i}(t)$ and $u_{i}(t)$ are adapted stochastic processes.

\begin{proposition}
Let $\mu$ be a relaxed control and $X$ the corresponding state process. Then
one can choose a sliding control%
\begin{equation}
\nu_{t}=%
{\displaystyle\sum\limits_{i=1}^{p}}
\alpha_{i}(t)\delta_{u_{i}(t)}(da),\text{ }u_{i}(t)\in A,\text{ }\alpha
_{i}(t)\geq0\text{ and }%
{\displaystyle\sum\limits_{i=1}^{p}}
\alpha_{i}(t)=1
\end{equation}

such that

1) $X$ is a solution of the controlled MFSDE%
\begin{equation}
\left\{
\begin{array}
[c]{lll}%
dX_{t} & = &
{\displaystyle\sum\limits_{i=1}^{p}}
\alpha_{i}(t)b(t,X_{t},E(\Psi(X_{t})),u_{i}(t))dt+%
{\displaystyle\sum\limits_{i=1}^{p}}
\alpha_{i}(t)^{1/2}\sigma(t,X_{t},E(\Phi(X_{t})),u_{i}(t))dW_{t}^{i}\\
X_{0} & = & x
\end{array}
\right.  \label{Sliding}%
\end{equation}

where $\left(  W^{i}\right)  $ are independant Brownian motions defined on an
extension of the probability space.

2) $J(\mu)=J(\nu),$ where $J(\nu)$ is associated with $X.$
\end{proposition}

\bop

Let $\Lambda$ denote the $d+d^{2}+1$-dimensional simplex

\begin{center}
$\Lambda=\left\{  \lambda=\left(  \lambda_{0},\lambda_{1},...,\lambda
_{d+d^{2}+1}\right)  ;\text{ }\lambda_{i}\geq0;\text{ }%
{\displaystyle\sum\limits_{i=0}^{d+d^{2}+1}}
\lambda_{i}=1\right\}  $
\end{center}

and $W$ the $(d+d^{2}+2)$-cartesian product of the set $\mathbb{A}$

\begin{center}
$W=\left\{  w=\left(  u_{0},u_{1},...,u_{d+d^{2}+1}\right)  ;\text{ }u_{i}%
\in\mathbb{A}\right\}  $
\end{center}

Define the function

\begin{center}
$g(t,\lambda,w)=%
{\displaystyle\sum\limits_{i=0}^{d+d^{2}+1}}
\lambda_{i}\widetilde{b}(t,X_{t},E(\Psi(X_{t})),u_{i})-%
{\displaystyle\int\nolimits_{A}}
\widetilde{b}\left(  t,X_{t},E(\Psi(X_{t})),a\right)  \mu_{t}(da)$
\end{center}

where $t\in\left[  0,T\right]  ,\lambda\in\Lambda,w\in W$ and $\widetilde
{b}(t,X_{t},E(\Psi(X_{t})),u_{i})=%
\begin{pmatrix}
b(t,X_{t},E(\Psi(X_{t})),u_{i})\\
\sigma\sigma^{\ast}(t,X_{t},E(\Phi(X_{t})),u_{i})\\
h(t,x_{t},E(\varphi(X_{t})),u_{i})
\end{pmatrix}
$ \bigskip

Let $\widetilde{b}(t,X_{t},E(\Psi(X_{t})),u_{i}),$ $i=0,1,...,d+d^{2}+1,$ be
arbitrary points in $P(t,X_{t})$ where

\begin{center}
$P(t,X_{t})=\left\{  \left(  b(t,X_{t},E(\Psi(X_{t})),a),\text{ }\sigma
\sigma^{\ast}(t,X_{t},E(\Phi(X_{t})),a),\text{ }h(t,X_{t},E(\Psi
(X_{t})),a)\right)  ;\text{ }a\in\mathbb{A}\right\}  \subset\mathbb{R}%
^{d+d^{2}+1}$
\end{center}

Then the convex hull of this set is the collection of all points of the form

\begin{center}
$%
{\displaystyle\sum\limits_{i=0}^{d+d^{2}+1}}
\lambda_{i}\widetilde{b}(t,X_{t},E(\Psi(X_{t})),u_{i})$
\end{center}

If $\mu$ is a relaxed control, then $%
{\displaystyle\int\nolimits_{\mathbb{A}}}
\widetilde{b}\left(  t,X_{t},E(\Psi(X_{t})),a\right)  \mu_{t}(da)\in
Conv\left(  P(t,X_{t})\right)  $, the convex hull of $P(t,X_{t}).$ Therefore
it follows from Carath\'{e}odory's Lemma (which says that the convex hull of a
$d$-dimensional set $M$ coincides with the union of the convex hulls of $d+1$
points of $M$), that for each $\left(  \omega,t\right)  \in\Omega\times\left[
0,T\right]  $ the equation $g(t,\lambda,\omega)=0$ admits at least one
solution. Moreover the set

\begin{center}
$\left\{  \left(  \omega,\lambda,w\right)  \in\Omega\times\Lambda\times W:%
{\displaystyle\sum\limits_{i=0}^{d+d^{2}+1}}
\lambda_{i}\widetilde{b}(t,X_{t},E(\Psi(X_{t})),u_{i})=%
{\displaystyle\int\nolimits_{\mathbb{A}}}
\widetilde{b}\left(  t,x_{t},E(\Psi(x_{t})),a\right)  \mu_{t}(da)\right\}  $
\end{center}

is measurable with respect to $\mathcal{F}_{t}\otimes\mathcal{B}%
(\mathbb{R}^{d+1})\otimes\mathcal{B}(\mathbb{A}^{d+1})$ with non empty
$\omega-$sections for each $\omega$.

Hence by using a measurable selection theorem \cite{EKNJ}, there exist
measurable $\mathcal{F}_{t}-$adapted processes $\lambda_{t}$ and $w_{t}$ with
values, respectively in $\Lambda$ and $W$ such that:

\begin{center}
$%
{\displaystyle\int\nolimits_{\mathbb{A}}}
\widetilde{b}\left(  t,X_{t},E(\Psi(X_{t})),a\right)  \mu_{t}(du)=%
{\displaystyle\sum\limits_{i=0}^{d+d^{2}+1}}
\lambda_{i}(t)\widetilde{b}(t,X_{t},u_{i}(t))$
\end{center}

Then it is easy to verify that the process defined by its drift $%
{\displaystyle\sum\limits_{i=0}^{d+d^{2}+1}}
\lambda_{i}(t)b(t,X_{t},E(\Psi(X_{t})),u_{i}(t))$ and its quadratic variation$%
{\displaystyle\sum\limits_{i=0}^{d+d^{2}+1}}
\lambda_{i}(t)\sigma\sigma^{\ast}(t,X_{t},E(\Psi(X_{t})),u_{i}(t))$ is the
solution of the MFSDE (\ref{Sliding}), defined possibly on an extension of the
initial probability space because of the possible degeneracy of the matrix
$\sigma\sigma^{\ast}$. \eop

\begin{corollary}
Assume that the set\textbf{ }%
\[
P(t,X_{t})=\left\{  \widetilde{b}(t,X_{t},E(\Psi(X_{t})),a);\text{ }%
a\in\mathbb{A}\right\}  \subset\mathbb{R}^{d+d^{2}+1}%
\]
is convex, where \bigskip$\widetilde{b}(t,X_{t},E(\Psi(X_{t})),a)=\left(
b(t,X_{t},E(\Psi(X_{t})),a),\text{ }\sigma\sigma^{\ast}(t,X_{t},E(\Phi
(X_{t})),a),\text{ }h(t,X_{t},E(\Psi(X_{t})),a)\right)  .$

Then the relaxed optimal control is realized by a strict control.
\end{corollary}

\bop

Using Proposition 2.19, it follows that for each relaxed control $\mu$ we have

\begin{center}
$%
{\displaystyle\int\nolimits_{\mathbb{A}}}
\widetilde{b}\left(  t,X_{t},E(\Psi(X_{t})),a\right)  \mu_{t}(da)\in
Conv\left(  P(t,X_{t})\right)  $
\end{center}

Since $P(t,X_{t})$ is convex then $Conv\left(  P(t,X_{t})\right)
=P(t,X_{t}).$ Then applying the same arguments as in Proposition 2.19, there
exists a measurable $\mathcal{F}_{t}-$adapted process $u_{t}$ such that

\begin{center}
$%
{\displaystyle\int\nolimits_{\mathbb{A}}}
\widetilde{b}\left(  t,X_{t},E(\Psi(X_{t})),a\right)  \mu_{t}(du)=\widetilde
{b}(t,X_{t},u_{t})$
\end{center}

which implies that $X_{t}$ is a solution of the MFSDE

\begin{center}
$\left\{
\begin{array}
[c]{l}%
dX_{t}=b(t,X_{t},E(\Psi(X_{t})),u_{t})dt+\sigma(t,X_{t},E(\Phi(X_{t}%
)),u_{t})dW_{t}\\
X_{0}=x
\end{array}
\right.  $
\end{center}

and $J(\mu)=J(u).$ This ends the proof. \eop

\section{Conclusion}

\textit{We have proved the existence of an optimal relaxed control for
mean-field stochastic control problems, where both the drift and diffusion
coefficient are controlled. The natural stochastic equation corresponding to a
relaxed control is a mean field stochastic equation driven by an orthogonal
martingale measure. As we have proved in a counter-example, replacing the
drift and the diffusion coefficient by their relaxed counterparts in the
relaxed equation, as in \cite{Cha}, does not lead to a true extension of the
original problem. In fact, the case where the diffusion coefficient is
controlled is not a direct extension of the deterministic case and reflects
the stochastic nature of the problem. }

\end{document}